\documentclass[11pt]{article}
\usepackage{amsmath}
\usepackage{mathrsfs}
\usepackage{amsfonts}

\usepackage{amsfonts,amsmath,amssymb}
\usepackage{amssymb,amsfonts,amsmath,
latexsym,epsfig,cite,psfrag,eepic,colordvi}
\usepackage{amscd,graphics}
\usepackage{lineno}

\newcommand{\qed}{\hfill $\Box $}
\newcommand{\pf}{\noindent {\bf Proof.} }

\newtheorem{theorem}{Theorem}[section]

\newtheorem{coro}[theorem]{Corollary}

\begin{document}

\title{On perfect $k$-matchings \thanks{This work is supported by the Fundamental Research Funds
for the Central Universities.}}

\author{Hongliang Lu\thanks{Corresponding email:
luhongliang215@sina.com (H. Lu)} 
\\ {\small Department of Mathematics}
\\ {\small Xi'an Jiaotong University, Xi'an 710049, PR China}
}

\date{}

\maketitle

\date{}

\maketitle

\begin{abstract}
In this paper, we generalize the notions of perfect matchings,
perfect 2-matchings  to perfect $k$-matchings and give a necessary
and sufficient condition for existence of perfect $k$-matchings. For
bipartite graphs, we show that this $k$-matching problem is
equivalent to that  matching question. Moreover, for regular graphs,
we  provide a sufficient condition
of perfect $k$-matching in  terms of edge connectivity.

\end{abstract}


\section{Introduction}

All graphs considered are multigraphs (with loops) and finite. Let
$G = (V(G),E(G))$ be a graph with vertex set $V(G)$ and edge set
$E(G)$. The number of vertices of a graph $G$ is called the
\emph{order} of $G$.  Unless otherwise defined, we follow
\cite{lovPlu} for terminologies and definitions.

 We denote the
degree of vertex $v$ in $G$ by  $d_{G}(v)$. For two subsets
$S,T\subseteq V(G)$, let $e_{G}(S,T)$ denote the number of edges of
$G$ joining $S$ to $T$. For a set $X$, we denote the cardinality of
$X$ by $|X|$. A vertex of degree zero is called an \emph{isolated
vertex}. Let $Iso(G)$ denote the  set of isolated vertices of $G$
and let $i(G)=Iso(G)$. Let $c_{o}(G)$ denote the number of odd
components of $G$. Let $odd(G)$ denote the number of odd components
with order at least three of $G$. For any subset $X$ of vertices of
$G$, we define the neighbourhood of $X$ in $G$ to be the set of all
vertices adjacent to vertices in $X$; this set is denoted by
$N_G(X)$.

A \emph{matching} $M$ of a graph $G$ is a subset of $E(G)$ such that
any two edges of $M$ have no end-vertices in common.
 Let
$k$ be a positive. A \emph{$k$-factor} of a graph $G$ is a spanning
subgraph $H$ of $G$  such that $d_{H}(x)=k$ for every $x\in V(G)$. A
\emph{$\{K_2,C_{2t+1}\ |\ t\geq 1\}$-factor} of a graph $G$ is a
spanning subgraph of $G$ such that each of its components is
isomorphic to one of $\{K_2,C_{2t+1}\ |\ t\geq 1\}$.


Let $f:\{0,1,\ldots,k\}\rightarrow E(G)$ be an assignment such that
the sum of weights of edges incident with any vertex is at most $k$,
i.e., $\sum_{e\sim v}f(e)\leq k$ for any vertex $v\in V(G)$.  A
$k$-matching is a subgraph induced by the edges with weight among
$1,\ldots,k$ such that $\sum_{e\sim v}f(e)\leq k$. The sum of all
weights, i.e., $\sum_{e\in E(G)}f(e)$, is called \emph{size} of a
$k$-matching $f$. A $k$-matching is \emph{perfect} if $\sum_{e\sim
v}f(e)=k$ for every vertex $v\in V(G)$. Clearly, a $k$-matching is
perfect if and only if its size is $k|V(G)|/2$. 
If $k=1$, then a perfect
$k$-matching is called a  \emph{perfect matching}. If $k=2$, then a
perfect $k$-matching is called a  \emph{perfect 2-matching}.

For perfect matching of  bipartite graphs, Hall  obtained the next
result in terms of isolated vertices.
\begin{theorem}[Hall, \cite{Hall}]\label{Hall}
Let $G = (X, Y )$ be a bipartite graph. Then $G$ has a perfect
matching if and only if $|X| = |Y |$ and for any $S\subseteq X$,
$$i(G-S)\leq |S|.$$
\end{theorem}

Tutte (1947) studied  the perfect matching of general graphs and
 gave the sufficient and necessary condition.

\begin{theorem}[Tutte, \cite{Tutte47}]\label{Tutte47}
A graph $G$ has a perfect matching if and only if for any
$S\subseteq V (G)$, $$ c_{o}(G-S)\leq |S|.$$
\end{theorem}

For perfect 2-matching, Tutte (1953) gave the following result.

\begin{theorem}[Tutte, \cite{Tutte53}]\label{Tutte53}
Let $G$ be a connected graph. Then the following statements are
equivalent:
\begin{enumerate}
\item [$(1)$] $G$ has a perfect 2-matching;

\item [$(2)$] $i(G-S)\leq |S|$ for all subsets $S\subseteq  V (G)$;

\item [$(3)$] $G$ has a $\{K_2,C_{2t+1}\ |\ t\geq 1\}$-factor.
\end{enumerate}
\end{theorem}

In the proof, we need the following technical theorems.

\begin{theorem}[Tutte, \cite{Tutte52}]\label{Tutte52}
 Let $G$ be a graph and $k$ a positive integer.
Then $G$ has a $k$-factor if and only if, for all $D, S \subseteq V
(G)$ with $D\cap S =\emptyset$,
\begin{align*}
\delta_{G}(D, S) = k|D|-k|S|+\sum_{v\in S}d_{G-D}(v)-\tau_G(S,
T)\geq 0,
\end{align*}
 where
$\tau_G(D, S)$ is the number of components $C$ of $G-(D\cup S$) such
that  $e_G(V (C), S)+ k|C|\equiv 1 \pmod 2$. Moreover, $\delta_G(D,
S)\equiv k|V (G)| \pmod 2$.

\end{theorem}

\section{Main Results}

In this section, we gave a good characterization for perfect
$k$-matchings.

\begin{theorem}\label{main1}
Let $k\geq 4$ be  even. Then $G$ contains a perfect $k$-matching if
and only if $G$ contains a perfect $2$-matching.
\end{theorem}

\pf Suppose that $G$ contains a perfect $2$-matching. By Theorem
\ref{Tutte53}, $G$ contains a $\{K_{2},C_{2l+1}\}$-factor $H$. We
assign every isolated edge of $H$ with weight $k$ and the rest edge
with weight $k/2$. Then we obtain a perfect $k$-matching of $G$.

Conversely, suppose $G$ that contains a perfect $k$-matching $H$.
Then there exists a function $f:V(G)\rightarrow \{0,1,\ldots,k\}$
such that $\sum_{v\sim e}f(e)=k$ for all $v\in V(G)$. We claim
$i(G-S)\leq |S|$ for all $S\subseteq V(G)$. Otherwise, assume that
there exists $S\subseteq V(G)$ such that $i(G-S)>|S|$.  Then we have
\begin{align*}
ki(G-S)=\sum_{e\in E_{G}(Iso(G-S),S)} f(e)>k|S|,
\end{align*}
a contradiction. 
So by Theorem \ref{Tutte53},
 $G$ contains a perfect $2$-matching. \qed

\begin{coro}
Let $k\geq2$ be  even. Then a graph $G$ contains a perfect
$k$-matching if and only if $i(G-S)\leq |S|$ for all $S\subseteq
V(G)$.
\end{coro}


\begin{theorem}\label{main2}
Let $k\geq 1$ be  odd. Then $G$ contains a perfect $k$-matching if
and only if
\begin{align*}
odd(G-S)+ki(G-S)\leq k|S|\ \  \ \ \mbox{for all subsets $S\subseteq
V(G)$}.
\end{align*}
\end{theorem}

\pf We first prove the necessity. Suppose that $G$ has a perfect
$k$-matching and there exists $S\subseteq V(G)$ such that
$$odd(G-S)+ki(G-S)>k|S|.$$
Let $f: E(G)\rightarrow \{0,1,\ldots,k\}$ such that $\sum_{e\sim
v}f(e)=k$ for all $v\in V(G)$.  Let $m=odd(G-S)$ and let
$C_{1},\ldots,C_{m}$ denote the odd components of $G-S$ with order
at least three. Let $W=C_{1}\cup\cdots \cup C_{m}$. Since $k$ is
odd, by parity, every odd component with order at least three can't
contain a perfect $k$-matching. So $\sum_{e\in
E_{G}(V(C_i),S)}f(e)\geq 1$ for $i=1,\ldots,m$.
 Then we have
\begin{align*}
k|S|=\sum_{v\in S}\sum_{e\sim v}f(e)&\geq \sum_{e\in
E_{G}(V(W),S)}f(e)+\sum_{e\in E_{G}(Iso(G-S),S)}f(e)\\
&\geq odd(G-S)+ki(G-S)>k|S|,
\end{align*}
a contradiction. So the result is followed.

We next prove the sufficiency.
 Let
$G^*$ be obtained from $G$ by changing every edge of $G$ into $k$
parallel edges. Then $G$ contains a perfect $k$-matching if and only
if $G^*$ contains a $k$-factor. Conversely, suppose that $G$
contains no  perfect $k$-matchings.  Then $G^*$ contains no
$k$-factors. By Theorem \ref{Tutte52}, there exist two disjoint
subset $D,S\subseteq V(G^*)$ such that
\begin{align*}
k|D|-k|S|+\sum_{x\in S}d_{G^*-D}(x)-\tau<0,
\end{align*}
where $\tau$ denote the number of components $C$ of $G^*-D-S$ such
that $k|V(C)|+e_{G^*}(V(C),S)\equiv 1$ (mod 2). Let
$C_{1},\ldots,C_{\tau}$ denote those components and
$W=\bigcup_{i=1}^\tau C_{i}$. By Theorem \ref{Tutte47}, we can
suppose that $k\geq 3$.

 Without loss of generality, among all
such subsets, we choose subsets $D$ and $S$ such that $S$ is
minimal. We have $S\neq\emptyset$, otherwise, $k|D|<\tau$ and
$|V(C_i)|$ is odd for $i=1\ldots,\tau$. So we have
\begin{align*}
k|D|-ki(G-D)-odd(G-D)\leq k|D| -\tau<0,
\end{align*} a contradiction.
Let $M=G^*-D-S-V(W)$.

\medskip
{\it Claim 1.} $G[S]$ consists of isolated vertices.
\medskip

Otherwise, let $e=uv\in G[S]$. Let $|N_G(v)\cap V(W)|=m$. Let $D'=D$
and $S'=S-v$. Let $\tau'$ denote the number of components $C$ of
$G^*-D'-S'$ such that $k|V(C)|+e_{G^*}(V(C),S')\equiv 1$ (mod 2).
Then we have
\begin{align*}
k|D'|-k|S'|+\sum_{x\in S'}d_{G^*-D'}(x)-\tau'&\leq
k|D|-k(|S|-1)+\sum_{x\in S-v}d_{G^*-D}(x)-(\tau-m)\\
&\leq k|D|-k|S|+k+\sum_{x\in S}d_{G^*-D}(x)-d_{G^*-D}(v)-(\tau-m)\\
&\leq k|D|-k|S|+k+\sum_{x\in S}d_{G^*-D}(x)-k(m+1)-\tau+m\\
&\leq k|D|-k|S|+\sum_{x\in S}d_{G^*-D}(x)-\tau<0,
\end{align*}
contradicting to the minimality of $S$. This completes the claim.

With the similar proof of Claim 1, we obtain the following claim.

\medskip
{\it Claim 2.} $e_{G}(S,V(M))=\emptyset$.
\medskip

\medskip
{\it Claim 3.} $|N_G(x)\cap V(W)|\leq 1$ for all $x\in S$.
\medskip

Otherwise, suppose that there exists $v\in S$ such that
$m=|N_G(v)\cap V(W)|\geq 2$. Let $D''=D$ and $S''=S-v$. Let $\tau''$
denote the number of components $C$ of $G^*-D''-S''$ such that
$k|V(C)|+e_{G^*}(V(C),S'')\equiv 1$ (mod 2). Then we have
\begin{align*}
k|D''|-k|S''|+\sum_{x\in S''}d_{G^*-D''}(x)-\tau''&\leq k|D|-k(|S|-1)+\sum_{x\in S-v}d_{G^*-D}(x)-(\tau-m)\\
&\leq k|D|-k|S|+k+\sum_{x\in S}d_{G^*-D}(x)-d_{G^*-D}(v)-(\tau-m)\\
&\leq k|D|-k|S|+k+\sum_{x\in S}d_{G^*-D}(x)-km-\tau+m\\
&= k|D|-k|S|+\sum_{x\in S}d_{G^*-D}(x)-\tau-(k-1)(m-1)+1\\
&\leq k|D|-k|S|+\sum_{x\in S}d_{G^*-D}(x)-\tau<0,
\end{align*}
contradicting to the minimality of $S$. This completes the claim.

\medskip
{\it Claim 4.} $E_{G}(S,V(W))=\emptyset$.
\medskip

Otherwise, by Claim 2, suppose that there exists an edge $uv\in
E_{G}(S,V(W))$, where $v\in S$ and $u\in V(W)$. Let $D'''=D$ and
$S'''=S-v$. Let $\tau'''$ denote the number of  components $C$ of
$G^*-D'''-S'''$ such that $k|V(C)|+e_{G^*}(S''',V(C))\equiv 1$ (mod
2). Without loss of generality, suppose that $u\in C_{1}$. By Claims
1, 2 and 3, then
  $G^*[V(C_{1})\cup \{v\}]$ is a component of $G^*-D'''-S'''$. Note that
   $k|V(C_{1})\cup \{v\}|+e_{G^*}(S''',V(C_{1})\cup \{v\})\equiv
 k|V(C_{1})|+e_{G^*}(S,V(C_{1}))\equiv 1$ (mod 2). 
 So $\tau=\tau'''$. Hence
\begin{align*}
k|D'''|-k|S'''|+\sum_{x\in S'''}d_{G^*-D'''}(x)-\tau'''
&= k|D|-k(|S|-1)+\sum_{x\in S-v}d_{G^*-D}(x)-\tau\\
&\leq k|D|-k|S|+k+\sum_{x\in S}d_{G^*-D}(x)-d_{G^*-D}(v)-\tau\\
&\leq k|D|-k|S|+\sum_{x\in S}d_{G^*-D}(x)-\tau<0,
\end{align*}
contradicting to the minimality of $S$. This completes the claim.

Since $e_{G^*}(V(C_i),S)+k|V(C_i)|\equiv 1\pmod2$ and $k$ is odd, by
Claim 4, we have $|V(C_i)|\equiv 1\pmod2$ for $i=1,\ldots,\tau$. By
Claims 1, 2, and 4, we have
\begin{align*}
0&>k|D|-k|S|+\sum_{x\in S}d_{G^*-D}(x)-\tau\\
&=k|D|-k|S|-\tau\\
&\geq k|D|-ki(G-D)-odd(G-D).
\end{align*}
Hence we have $k|D|<ki(G-D)+odd(G-D)$, a contradiction. We complete
the proof. \qed

%
%

\begin{theorem}
Let $G=(U,W)$ be a bipartite graph, where $|U|=|W|$. Then $G$
contains a perfect matching if and only if $G$ contains a perfect
$k$-matching.
\end{theorem}

\pf  Necessity is obvious. Now we prove the sufficiency.    
Suppose that $G$ contains a perfect $k$-matching. Let
$f:E(G)\rightarrow \{0,1,\ldots,k\}$ such that $\sum_{v\sim
e}f(e)=k$ for all $v\in V(G)$. Then for all independent set $S$, we
have
\begin{align}
k|S|=\sum_{v\in S}\sum_{v\sim e}f(e)\leq \sum_{ e\in
E_{G}(S,N(S))}f(e)\leq k|N(S)|.
\end{align}
So we have $i(G-S)\leq |S|$ for all $S\subseteq U$. By Theorem
\ref{Hall}, $G$ contains a perfect matching. This completes the
proof. \qed

\begin{coro}
Let $k\geq 1$ be an odd integer and $G$ be an $r$-regular,
$\lambda$-edge-connected  graph. Suppose that
$$
\lambda = \begin{cases} \lceil \frac{r}{k}\rceil-1 \ \ \
\qquad\mbox{ if } \lceil \frac{r}{k}\rceil\equiv r \pmod 2;\\
 \lceil \frac{r}{k}\rceil  \qquad \qquad\ \
\mbox{ if }  \lceil \frac{r}{k}\rceil\neq r \pmod 2.\end{cases}
$$
Then $G$ contains a perfect  $k$-matching.
\end{coro}

\pf  Suppose that the result doesn't hold.
 By Theorem
\ref{main2}, there exists a subset $S\subseteq V(G)$ such that
\begin{align*}
odd(G-S)+ki(G-S)>k|S|.
\end{align*}
 Let $m=odd(G-S)$. Let $C_1,\ldots,C_m$ denote
these odd components  with order at least three of $G-S$. Since
$r|C_i|-e_{G}(V(C_i),S)=\sum_{x\in V(C_i)}d_{C_i}(x)$ is even, so
$e_{G}(V(C_{i}),S)\equiv r$. Since $G$ is an $r$-regular,
$\lambda$-edge-connected graph, so
 if $\lceil \frac{r}{k}\rceil\equiv r
\pmod 2 $, then $e_{G}(V(C_{i}),S)\geq \lambda+1$ for
$i=1,\ldots,m$. So we have
\begin{align*}
r|S|\geq \lceil \frac{r}{k}\rceil odd(G-S)+ri(G-S).
\end{align*}
Hence,
\begin{align*}
kr|S|&\geq k\lceil \frac{r}{k}\rceil odd(G-S)+kri(G-S)\\
&\geq r odd(G-S)+k r i(G-S)>kr|S|,
\end{align*}
a contradiction. This completes the proof. \qed

\begin{coro}[B\"abler, \cite{Babler}]
Let  $G$ be an $r$-regular, $(r-1)$-edge-connected  graph with even
order. Then $G$ contains a perfect  matching.
\end{coro}
%




\end{document}